\documentclass[11pt]{amsart}
\usepackage{amssymb,amscd,amsmath}
\usepackage[all]{xy}
\newtheorem{thm}{Theorem}[section]

\newtheorem{cor}[thm]{Corollary}

\newtheorem{prop}[thm]{Proposition}
\newtheorem{lemma}[thm]{Lemma}

\newtheorem{conj}[thm]{Conjecture}
\newtheorem{prob}[thm]{Problem}

\newtheorem{lem}[thm]{Lemma}

\theoremstyle{remark}

\theoremstyle{definition}

\numberwithin{equation}{section}

\renewcommand{\bar}{\overline}

\newcommand{\cG}{\mathcal{G}}

\newcommand{\cS}{\mathcal{S}}
\newcommand{\cU}{\mathcal{U}}

\newcommand{\cX}{\mathcal{X}}
\newcommand{\cY}{\mathcal{Y}}

\newcommand{\Spin}{\mathrm{Spin}}
\newcommand{\SO}{\mathrm{SO}}

\newcommand{\compose}{{\scriptstyle\circ}}

\newcommand{\A}{{\mathbb{A}}}

\renewcommand{\P}{{\mathbb{P}}}
\newcommand{\F}{{\mathbb{F}}}

\newcommand{\Z}{{\mathbb{Z}}}
\newcommand{\tr}{\mathrm{tr}}

\newcommand{\fix}{\mathrm{fix}}
\newcommand{\cyc}{\mathrm{cyc}}

\newcommand{\SL}{\mathrm{SL}}
\newcommand{\Sp}{\mathrm{Sp}}

\newcommand{\Span}{\mathrm{Span}}

\newcommand{\Spec}{\mathrm{Spec}\;}
\newcommand{\id}{\mathrm{id}}

\newcommand{\ad}{\mathrm{ad}}

\newcommand{\g}{\mathfrak g}

\newcommand{\m}{\mathfrak m}
\newcommand{\p}{\mathfrak p}

\newcommand{\Alt}{{\raise 2pt\hbox{$\scriptstyle\bigwedge$}}}

\newcommand{\la}{\lambda}

\begin{document}
\title[Word maps]
{Word maps and Waring type problems}

\author{Michael Larsen}
\email{larsen@math.indiana.edu}
\address{Department of Mathematics\\
    Indiana University \\
    Bloomington, IN 47405\\
    U.S.A.}

\author{Aner Shalev}
\email{shalev@math.huji.ac.il}
\address{Einstein Institute of Mathematics\\
    Hebrew University \\
    Givat Ram, Jerusalem 91904\\
    Israel}

\thanks{Michael Larsen was partially supported by NSF grant DMS-0354772.
Aner Shalev was partially supported by an Israel Science Foundation
Grant. Both authors were partially supported by a Bi-National Science 
Foundation United States-Israel Grant.}

\begin{abstract}
Waring's classical problem deals with expressing every natural 
number as a sum of $g(k)$ $k$th powers. Recently there has been
considerable interest in similar questions for nonabelian groups,
and simple groups in particular. Here the $k$th power word
is replaced by an arbitrary group word $w \ne 1$, and the goal
is to express group elements as short products of values of
$w$.

We give a best possible and somewhat surprising solution for 
this Waring type problem
for various finite simple groups, showing that a product of
length two suffices to express all elements. We also show that
the set of values of $w$ is very large, improving various
results obtained so far.

Along the way we also obtain new results of independent interest 
on character values and class squares in symmetric groups.

Our methods involve algebraic geometry, representation theory,
probabilistic arguments, as well as three prime theorems from 
additive number theory (approximating Goldbach's Conjecture).

\bigskip

Contents:

1. Introduction

2. Word values in groups of Lie type

3. Waring's problem in bounded dimension

4. Special word values in $\SL_2(\F_p)$

5. Alternating groups, I: three primes theorems

6. Alternating groups, II: a combinatorial construction

7. Alternating groups, III: characters and probability 

8. Intersection theorems

\end{abstract}

\maketitle

\newpage

\section{Introduction}

Waring's problem asks whether every natural number is a sum 
of $g(k)$ $k$th powers (where $g$ is a suitable function). 
This was solved affirmatively by Hilbert in 1909.  
Optimizing $g(k)$ has been a 
central problem in additive number theory ever since 
(see \cite{Na} for more detail and background).

Recently there has been considerable interest in group theoretic analogues
of this phenomenon, where the aim is to present group
elements as short products of certain ``special'' elements.
These special elements can be powers, or commutators, or
values of a general word $w$, or elements of a special
conjugacy class in the group.

The motivation for this is sometimes topological.
Let $G$ be a (topologically) finitely generated pro-$p$ group.
Using expressions of elements of the commutator subgroup $G'$ 
as bounded products of commutators Serre showed that $G'$
is closed, and deduced that every finite index subgroup of 
$G$ is open. A recent deep result of Nikolov and Segal 
\cite{NS1, NS2} shows that this holds for every finitely 
generated profinite group. Again the core of the proof is
presenting group elements as short products of values
of certain words.

In the realm of finite simple groups such Waring type problems
have been studied in even greater detail. We need some notation.

Let $w = w(x_1, \ldots , x_d)$ be a non-trivial 
group word, namely a non-identity element of the free group $F_d$ on 
$x_1, \ldots, x_d$. Then we may write
$w = x_{i_1}^{n_1} x_{i_2}^{n_2} \cdots x_{i_k}^{n_k}$
where $i_j \in \{ 1, \ldots , d \}$, $n_j$ are integers, and we may
assume further that $w$ is reduced.
Let $G$ be a group. For $g_1, \ldots , g_d \in G$ we write
\[
w(g_1, \ldots , g_d) = g_{i_1}^{n_1} g_{i_2}^{n_2} \cdots g_{i_k}^{n_k} 
\in G.
\]
Let
\[
w(G) = \{ w(g_1, \ldots , g_d): g_1, \ldots , g_d \in G \}
\]
be the set of values of $w$ in $G$. 
For subsets $A, B \subseteq G$ let $AB = \{ ab | a \in A, b \in B \}$ 
and $A^k = \{ a_1 \cdots a_k | a_i \in G \}$.

Fix a non-trivial group word $w$ and let $G$ be a finite simple
group. If $G$ is large enough then it follows from Jones \cite{J}
that $w(G) \ne \{ 1 \}$ (namely $w$ is not an identity in $G$).
Can we then find a constant $c$ (which may depend on $w$ but not
on $G$) such that $w(G)^c = G$?

Various instances of this problem have been solved affirmatively 
in the past decade or two. See Wilson \cite{W} for the commutator word 
$w(x_1,x_2) = [x_1,x_2] = x_1^{-1}x_2^{-1}x_1x_2$, and 
Martinez and Zelmanov \cite{MZ} and Saxl and Wilson \cite{SW} 
for power words $w= x_1^k$.
It follows from their result that every element of a large enough finite 
simple group is a product of $f(k)$ $k$th powers.

In \cite{LiSh} this is generalized to arbitrary words $w$.
It is proved there that if $w(G) \ne \{ 1 \}$ then $w(G)^c = G$, where
the exponent $c = c(w)$ depends on $w$ (but is not given
explicitly). 
A somewhat surprising new result from \cite{Sh} is that
the exponent $c$ can be chosen independent of the word $w$ and may be, in fact,
a very small number.
More specifically {\it for every word $w \ne 1$ there is
a positive integer $N = N(w)$ such that for every finite
simple group $G$ with $|G| \ge N(w)$ we have} 
\[
w(G)^3 = G.
\]

The main purpose of this paper is to prove an even stronger
result for certain families of finite simple groups, where 
the exponent $3$ is replaced by $2$ (see Theorems \ref{two-words}
and \ref{bounded-waring} below).
Along the way we prove several results of independent interest,
related to the size of $w(G)$, and to powers of certain conjugacy
classes $C$ in $G$.

Given a conjugacy class $C \ne 1$ in a finite simple group $G$
there exists a number $k$ such that $C^k=G$. Substantial work
has been devoted to the study of these numbers $k$, and related
so-called covering numbers (see e.g. \cite{AH, EGH, LL, LiSh}).

A particular challenge is to show that $C^2= G$ in certain cases.
Indeed a conjecture of Thompson states that every finite simple group
$G$ has a class $C$ with this property. In spite of considerable
progress this is still open in general.

Now let $C = \sigma^{S_n}$ be a conjugacy class in $S_n$. 
When can we say that $C^2 = A_n$? 
This problem has quite a long history.
Gleason \cite{Hu} seems to have been the first to observe that
$C^2 = A_n$ if $\sigma$ is an $n$-cycle.
See also Bertram \cite{Be}.
Later this was generalized by Brenner \cite{Br} to the case where $\sigma$ consists
of two cycles (or more generally two non-trivial cycles with some additional fixed points).
The case of permutations with more general cycle structure remained 
wide open.

Our first result deals with permutations with any given
number of cycles, provided $n$ is sufficiently large.
We denote the number of cycles (including 1-cycles) of a permutation 
$\sigma$ by $\cyc(\sigma)$.

\begin{thm} 
\label{squares}
For every positive integer $k$ there is a number $f(k)$ such that 
if $n \ge f(k)$ and $\sigma$ is a permutation in $S_n$ with 
$\cyc(\sigma) \le k$ then $(\sigma^{S_n})^2 = A_n$.

Moreover, there is an absolute constant $c$ such that,
if $n \ge c$ and $\sigma$ is a permutation in $S_n$
with $\cyc(\sigma) \le n^{1/128}$,
then $(\sigma^{S_n})^2 = A_n$.
\end{thm}

By the Erd{\H o}s-Tur{\'a}n theory (see \cite{ET}), a randomly 
chosen permutation in $S_n$ has $O(\log n)$ cycles. 
This gives rise to the following result.

\begin{cor} 
\label{randomconj}

(i) The probability that a randomly chosen permutation
$\sigma \in S_n$ satisfies $(\sigma^{S_n})^2 = A_n$ tends to
$1$ as $n \rightarrow \infty$. 

(ii) The probability that a randomly chosen permutation
$\sigma \in A_n$ satisfies $(\sigma^{A_n})^2 = A_n$ tends to
$1$ as $n \rightarrow \infty$. 

\end{cor}

Indeed, part (i) follows immediately, and part (ii) follows from (i)
using the fact that for almost all $\sigma \in A_n$ we have
$\sigma^{A_n} = \sigma^{S_n}$ (see \cite{Sh}, 3.3).

Corollary \ref{randomconj} sharpens a previous result showing that
for almost all $\sigma \in A_n$ we have $(\sigma^{A_n})^3 = A_n$ 
(see Theorem 2.6 in \cite{Sh}).

It is intriguing that the proof of Theorem~\ref{squares} is non-elementary
in the sense that it involves probabilistic arguments and new
character-theoretic estimates. In particular we show that all character
values $\chi(\sigma)$ of a permutation $\sigma \in S_n$ can be bounded
in terms of $\cyc(\sigma)$ alone (see Theorem~\ref{char-bound}).

In order to apply Theorem~\ref{squares} to Waring type problems 
we need to show that the image $w(A_n)$ of a word map $w \ne 1$
contains conjugacy classes $C$ with few cycles.

\begin{thm} 
\label{cyc}
For every positive integer $n$ there exists a permutation
$\sigma_n \in A_n$ with the following properties:

(i) $\cyc(\sigma_n)\le 23$;

(ii) For each non-trivial word $w$ there exists $N = N(w)$ such
that if $n \ge N$ then $\sigma_n \in w(A_n)$.
\end{thm}

We can also arrange that, for $n$ sufficiently large, the
permutation $\sigma_n$ has at least six fixed points.

Our proof of Theorem \ref{cyc} 
is highly non-elementary, involving algebraic geometry
and results from analytic number theory (such as weak versions
of the Goldbach Conjecture). The idea is to embed groups
of the form $\SL_2(\F_p)$ and their products in $A_n$,
and then focus on the properties of word maps on $\SL_2(\F_p)$.

It would be interesting to try to find an elementary proof
of Theorem \ref{cyc}, as well as to improve it.
In this context we pose the following:

\begin{prob} What is the minimal number $k$ such that
(for $w \ne 1$ and $n \ge N(w)$) $w(A_n)$ always contains 
permutations with at most $k$ cycles?
Does $w(A_n)$ necessarily contain a cycle of length
$n-k$ for some $k$ which may depend on $w$ but not on $n$?
\end{prob}

Equipped with the tools above, we can now improve earlier Waring
type results for alternating groups, showing that $w(A_n)^2 = A_n$
if $n \gg 0$.
In fact we prove a bit more.

\begin{thm}
\label{two-words}
For each pair of non-trivial words $w_1, w_2$ there exists 
$N = N(w_1, w_2)$ such that for all integers $n \ge N$ we have
\[
w_1(A_n)w_2(A_n) = A_n.
\]
\end{thm}

The idea of the proof is to combine Theorems \ref{squares}
and \ref{cyc}, noting that $\sigma_n^{S_n} \subset w_i(A_n)$
for $i = 1,2$ and $(\sigma_n^{S_n})^2 = A_n$.

We next establish a similar result for certain finite simple
groups of Lie type. 
Indeed, if we limit our attention to groups of Lie type of bounded 
dimension, we can prove the following result:

\begin{thm}
\label{bounded-waring}
Given an integer $d$ and two non-trivial words $w_1$ and
$w_2$ there exists an integer $N=N(d,w_1,w_2)$
such that if $G$ is a simply connected almost simple algebraic group
of dimension $d$ over a finite field $\F$, 
$\Gamma = G(\F)/Z(G(\F))$ is the finite simple group associated to 
$G$ over $\F$, and $|\Gamma|\ge N$, then
we have
$$w_1(\Gamma)w_2(\Gamma) = \Gamma.$$
\end{thm}

The method of the proof is to obtain a suitable surjectivity theorem 
at the level of algebraic groups,
and then to use the Riemann hypothesis for varieties over finite fields to show 
that the relevant fibers actually have points over $\F$.   

As a consequence we obtain the following.

\begin{cor} 
\label{three-words}
For each triple of non-trivial words $w_1, w_2, w_3$ there exists 
$N = N(w_1, w_2, w_3)$ such that if $G$ is a finite simple group
of order at least $N$, then
\[
w_1(G)w_2(G)w_3(G) = G.
\]
\end{cor}

Indeed, the case of groups of Lie type appears in \cite{Sh},
and Theorem \ref{two-words} above completes the proof by dealing 
with alternating groups.

We conjecture this can be strengthened as follows.

\begin{conj}
\label{challenge}
For each pair of non-trivial words $w_1, w_2$ there exists 
$N = N(w_1, w_2)$ such that if $G$ is a finite simple group
of order at least $N$, then
\[
w_1(G)w_2(G) = G.
\]
\end{conj}

In view of our results above it remains to prove this
for classical groups of unbounded rank. This seems quite
challenging even in the particular case $w_1 = w_2 = x_1^2$.

Let us now discuss the size of the subset $w(G)$.
A related useful result of Borel \cite{Bo} shows that the word map
induced by $w$ on simple algebraic groups is a dominant map.
In \cite{L} this is used to show that if $G$ is a finite simple group
and $\epsilon > 0$ then 
\[
|w(G)| \ge |G|^{1-\epsilon}
\]
provided $|G| \ge f(w,\epsilon)$.
In this paper we improve this bound
for various families of finite simple groups.

We start with alternating groups.
Often the key to estimating the size of $w(G)$ is finding
special elements in it. 

Indeed, using Theorem~\ref{cyc} and the remark following it 
one can easily deduce that, for $n \ge N(w)$,
$$|w(A_n)|\ge c n^{-17}|A_n|$$
for some absolute constant $c > 0$.  
In fact, we can provide a much better estimate:

\begin{thm} 
\label{lower-bound}
For each non-trivial word $w$ and every $\epsilon > 0$,
there exists $N = N(w, \epsilon)$ such that if $n \ge N$ then
\[
|w(A_n)| \ge n^{-4-\epsilon}|A_n|.
\]
\end{thm}

We can show that this bound is tight up to the value of the
exponent. Indeed, consider a power word $w(x_1) = x_1^k$.
Then by \cite{L2} we have
\[
|w(A_n)| \le |w(S_n)| \le a n^{-b} n!,
\]
where $b = 1 - \phi(k)/k$ and $a$ is some constant (depending
on $k$). Note that, choosing $k$ suitably (e.g. as the product
of the first $m$ primes) we may arrange that $b \ge 1-\epsilon$
for any fixed $\epsilon > 0$. This shows that the exponent 
in Theorem~\ref{lower-bound} must be at most $-1$.  

It would be interesting to find out whether 
Theorem~\ref{lower-bound} can be improved as follows.

\begin{prob}
\label{strong} Is it true that $|w(A_n)|\ge n^{-1} |A_n|$
provided $w \ne 1$ and $n \ge N(w)$?
\end{prob}

While we are unable to settle this for alternating groups,
we obtain an analogous result for certain simple groups of
Lie type.

\begin{thm} 
\label{not-A}
Let $G$ be a finite simple group
of Lie type and of rank $r$. Let $w \ne 1$ be a word. Then
if $G$ is not of type $A_r$ or $^2A_r$, we have
\[
|w(G)| \ge cr^{-1}|G|
\]
for some absolute constant $c>0$, provided $|G| \ge N(w)$.
\end{thm}

It is interesting that our methods allow us to generalize
various results on $w(G)$ to intersections of the form 
$\cap_{i=1}^k w_i(G)$, where $w_1, \ldots , w_k \ne 1$ are any
given words.
These generalizations are formulated in the last section of this
paper.

Some words on the structure of this paper.
Section 2 is devoted to groups of Lie type
and the proof of Theorem~\ref{not-A}.
In Section 3 we focus on Lie type groups of bounded dimension
and prove Theorem~\ref{bounded-waring}.
In Section 4 we make a closer analysis of groups of the
type $\SL_2(\F_p)$. Combined with a certain three primes theorem
this leads to the proofs of Theorem~\ref{cyc} and Theorem~\ref{lower-bound} 
in Section 5.
Section 6 contains the elementary part of the proof of Theorem 
\ref{squares}.
In Section 7 we study character values of $S_n$ on permutations
with few cycles and squares of conjugacy classes of such
permutations.  In this way we obtain a probabilistic proof of 
Theorem~\ref{squares}. 
Finally, in Section 8 we formulate results on the intersection
of the images of different word maps, and show how to derive them
using our methods.

The first-named author would like to thank Nicholas Katz for 
some useful correspondence.

\bigskip

\section{Word values in groups of Lie type}

In this section, we use algebraic geometry and group-theoretic arguments
to give lower bounds on $|w(\Gamma)|$ for certain groups $\Gamma$ of
Lie type, thus proving Theorem~\ref{not-A}.
Note that the bounds obtained are actually stronger than those we give 
for alternating groups.   
Although in the rest of the paper, $G$ denotes a finite simple group, in this
section and \S3, $G$ will always
denote a simply connected, almost simple algebraic group
defined over a finite field.  To get a finite simple group $\Gamma$,
we divide $G(\F_q)$ by its center.  All simple groups of Lie type arise in this way, excepting
the Suzuki and Ree groups.  The Suzuki and Ree groups can also be treated by the methods
of this section, but we exclude them to avoid unpleasant technicalities.

\begin{prop}
\label{fixed-rank}
Let $w$ be a non-trivial word and $G$ a simply connected almost simple 
algebraic group
over a finite field $\F_q$.  Then there is a positive constant $c$ depending
only on $w$ and $\dim G$ such that 
$$|w(G(\F_q))| > c|G(\F_q)|.$$
\end{prop}

This proposition is proved as \cite[Prop. 7]{L} in a form that also
includes the case of Ree and Suzuki groups.  Here we give an alternative
proof that uses standard techniques from algebraic geometry instead of
results from \cite{LP}.

We use the following more or less standard result:

\begin{lemma}
\label{alg-geom}
Let $\pi\colon \cX\to\cY$ denote a dominant morphism between reduced schemes of finite type over
$\Z$.  Then there exist constants $c_1$ and $c_2$ such that for every prime $p$ and every
point $y \in \cY(\F_p)$, we have
$$|\{x\in \cX(\F_p)\colon \pi(x) = y\}| \le c_1 p^{\dim \pi^{-1} x}$$
and if $\pi^{-1} y$ is non-empty and geometrically irreducible, then
$$|\{x\in \cX(\F_p)\colon \pi(x) = y\}| 
\ge p^{\dim \pi^{-1} y}(1-c_2 p^{-1/2}).$$
\end{lemma}

\begin{proof}
By the Lefschetz trace formula, these results are immediate consequences
of standard facts about \'etale cohomology groups:
the boundedness of the compactly supported cohomology groups
of the fibers of $\pi$, the trace map on the top-dimensional cohomology of
(geometrically irreducible) varieties, and the weight estimates
for the eigenvalues of Frobenius on these cohomology groups
(see \cite[XVIII~2.9]{SGA4} and \cite[3.3.1]{De}.)
\end{proof}

We now prove the proposition.

\begin{proof}
An algebraic group $G$ of this kind is specified by a connected 
Dynkin diagram $\Delta$,
an automorphism of this diagram (which is determined up to
inner automorphism by its order $g\le 3$), and a finite field $\F_q$.
It suffices to prove this proposition for each fixed choice of $\Delta$ and $g$.

First we treat the split case, $g=1$.
In this case, $\Delta$ determines a simply connected
Chevalley scheme $\cG$ over $\Spec \Z$ of which $G$ is
the $\F_q$-fiber, and
$w$ determines a morphism of schemes $\cG^d\to \cG$ which 
we denote $\pi_w$.
If $\pi_w$ has all fibers of dimension $\le(d-1)\dim G$, then we are done by
Lemma~\ref{alg-geom}, since 
$|G(\F_q)|\ge c_1 q^{\dim G}$ 
and each fiber has at most $c_2 q^{(d-1)\dim G}$ points of
$G(\F_q)^d$ for some positive constants $c_1$ and $c_2$ independent
of $q$.  The difficulty is that it might happen, a priori, that
most points in $G(\F_q)^d$ lie 
in fibers of $\pi_w$ of dimension greater than $(d-1)\dim G$.

To show that this does not in fact
happen for large $q$ we note first that $\pi_w$
is dominant; moreover, it is dominant for every fiber of $\cG^d$
over $\Spec\Z$ \cite{Bo}.  Fiber dimension
is a constructible function \cite[9.5.5]{EGA4-3}, 
so the Zariski closure of the set of geometric points $\bar s$
for which 
$$\dim \pi_{w}^{-1}(\bar s) > (d-1)\dim G$$
is a proper 
closed subset of $\cG$.  Let $\cX$ denote its inverse image in $\cG^d$.
By \cite{Bo}, $\cX$ is a proper closed subvariety of
$\cG^d$ for all $p$ which does not contain the generic point of
any fiber of $\cG^d$ over $\Spec\Z$.
By Lemma~\ref{alg-geom}, for all $q\gg 0$,
at least half of the points of $G(\F_q)^d$ do not lie in 
$\cX(\F_q)$, and the proposition follows.

Finally, we treat the case $g\ge 2$.  
Instead of defining a split group scheme $\cG$ over
$\Spec \Z$, we define a non-split group scheme $\pi\colon\cG\to\cS$ 
for some higher-dimensional base $\cS$
such that $\cG$ splits over a finite \'etale Galois 
cover $\cS'$ over $\cS$ of degree $g$.
For example, we may take 
$$\cS=\Spec \Z[1/g][x^g,x^{-g}]\coprod \Spec \F_g[x^g-x]$$
and 
$$\cS' = \Spec\Z[1/g][x,x^{-1}]\coprod \Spec \F_g[x].$$
The crucial point is that the base
is large enough that for every $\F_q$ there exists an $\F_q$-valued point $s$ of $\cS$ whose preimage
in $\cS'$ is $\F_{q^g}$, 
since this guarantees that every $G/\F_q$ with the given Dynkin diagram and diagram automorphism
can be obtained as the fiber of $\cG$ for some $\F_q$ point of $\cS$.

\end{proof}

We can now prove Theorem~\ref{not-A}.

\begin{proof}
The case of groups of type E, F, and G and of classical groups of any
fixed dimension follows from Proposition \ref{fixed-rank}.  
We need therefore consider only orthogonal groups 
and symplectic groups.  We begin with the symplectic case.

Let $V$ be a $2$-dimensional vector space over the field $\F_{q^r}$.
We endow $V$ with an area form $A$.  We define a bilinear form on $V$,
regarded as a vector space over $\F_q$, as follows:
$$\langle v_1,v_2\rangle = \tr_{\F_{q^r}/\F_q} A(v_1\wedge v_2).$$
This form is obviously antisymmetric and non-degenerate.  It therefore 
determines an inclusion 
$$i\colon \SL_2(\F_{q^r}) \to \Sp_{2r}(\F_q).$$
If
$x\in \SL_2(\F_{q^r})$ has eigenvalues $\lambda^{\pm 1}$,
then $i(x)$ has eigenvalues $\lambda^{\pm q^k}$,
$k=0,1,2,\ldots,r-1$.
As long as these eigenvalues are pairwise distinct, $i(x)$
is a regular semisimple symplectic matrix.  This can be achieved by insisting
for every proper subfield $\F_{q^d}$, $\lambda$ belongs neither
to $\F_{q^d}^\times$ nor to 
the norm-$1$ subgroup of $\F_{q^{2d}}$ over $\F_{q^d}$.  
The subset of $\SL_2(\F_{q^r})$ consisting of elements which
violate either condition has cardinality $O(q^{3r/2})$, so by
Proposition~\ref{fixed-rank}, we can find at least $c_1 q^{3r}$ elements
of $w(\SL_2(\F_{q^r}))$ which are regular semisimple in $\Sp_{2r}(\F_q)$.
If semisimple elements $a,b\in \SL_2(\F_{q^r})$ have eigenvalues 
$\alpha^{\pm1}$, $\beta^{\pm 1}$ respectively, 
then $i(a)$ and $i(b)$ have the same eigenvalues if and only if 
$\alpha = \beta^{q^i}$ or $\alpha^{-1} = \beta^{q^i}$ for some
$i\in\{0,1,\ldots,r-1\}$.  Thus there are at least $c_2 r^{-1} q^r$
distinct conjugacy classes of regular semisimple elements in
$i(w(\SL_2(\F_{q^r})))\subset w(\Sp_{2r}(\F_q))$.  The centralizer of
any regular semisimple element $i(a)$ in $\Sp_{2r}(\F_q)$ coincides with the
image of the centralizer of $a$ in $\SL_2(\F_{q^r})$; it is therefore
of order $q^r\pm 1$.  Thus, 
$$|i(w(\SL_2(\F_{q^r})))| \ge c_2 r^{-1} \frac{q^r}{q^r+1}|\Sp_{2r}(\F_q)|.$$
This proves the theorem for groups of type C.

For the cases B and D, we begin by observing that given 
a vector space $V$ over $\F_{q^k}$
and a non-degenerate quadratic form $Q_0\colon V\to \F_{q^k}$ defined over
$\F_{q^k}$,
we can define a non-degenerate quadratic form on
$V$ regarded as an $\F_q$-vector space as follows:
$$Q(v) = \tr_{\F_{q^k}/\F_q}Q_0(v).$$
This gives an inclusion $\Spin(V,Q_0)<\Spin(V,Q)$.
Composing this map with an isomorphism 
$\SL_2(\F_{q^{2k}})\cong \Spin_4^-(\F_{q^k})$,
we obtain an inclusion 
\begin{equation}
\label{o4n}
\SL_2(\F_{q^{2k}})\to \Spin_{4k}^{\pm}(\F_q)
\end{equation}
for some choice of signs.  
By the classification of quadratic forms over
finite fields, every space $W$ of dimension $\ge 4$ with a non-degenerate
quadratic form over $\F_q$ can be decomposed as
$W_1\perp W_2$, where $W_1$ is isomorphic to a space $V$ obtained as above and
$\dim W_2\le 6$.  We choose such a decomposition, and let $i$ denote the
map 
$$\SL_2(\F_{q^{2k}})\times \Spin(W_2)\to \Spin(W)$$ 
obtained by combining
(\ref{o4n}) with the natural map $\Spin(W_1)\times\Spin(W_2)\to \Spin(W)$.
Note that $i$ is at most two-to-one since $\SO(W_1)\times \SO(W_2)\to \SO(W)$
is injective.
If $k\ge 4$, $b\in \Spin(W_2)$ is regular semisimple, 
and the eigenvalues of a semisimple element 
$a\in \SL_2(\F_{q^{2k}})$ are not contained in a proper subfield
of $\F_{q^{2k}}$, then the image of $i(a,b)$ in $\SO(W)$
has all eigenvalues distinct and is therefore
regular semisimple.  By Proposition~\ref{fixed-rank}, the number of such pairs 
$(a,b)$ with $a\in w(\SL_2(\F_{q^{2k}}))$ is at least
$$c_3 q^{2k} q^{\frac{\dim W_2}2} = c_3 q^{\frac{\dim W}2},$$
and the number of semisimple conjugacy classes in $w(\Spin(W))$ determined
by such pairs is at least $c_4 r^{-1} q^{\frac{\dim W}2}$.  
For $q$ greater than some constant 
$c_5$ (independent of $k$), we can further specify that $b\in w(\Spin(W_2))$ at the cost
of replacing $c_4$ by $c_6$.  The centralizer of $i(a,b)$ in $\Spin(W)$ is 
a maximal torus of
$\Spin(W)$ containing a maximal torus of $\Spin(W_1)$; it therefore has at most 
$$(q^{2k}+1)(q+1)^{\frac{\dim W_2}2} \le c_7 q^{\frac{\dim W}2}$$
elements.  We conclude that the total number of elements 
of $\Spin(W)$ conjugate to elements of 
$i(\SL_2(\F_{q^{2k}})\times \Spin(W_2))$ is at least $\frac {c_6}{c_7}r^{-1}|\Spin(W)|$.

Finally, we consider the case that $q\le c_5$.  In this case, we take $b=1$, so $i(a,b)$
is not regular; its centralizer is $\Spin(W_2)$ times a maximal torus of $\Spin(W_1)$.
The order of the centralizer is therefore less than $c_5^{(\dim W_2)^2}q^{\frac{\dim W}2}$,
and the argument goes through as before.

\end{proof}

\bigskip

\section{Waring's problem in bounded dimension}

In this section we focus on simple groups of Lie type $G$ 
in bounded dimension. Using geometric methods 
and establishing the irreducibility of certain fibers (see 
Theorem~\ref{basic-irred} below), we obtain a best possible 
solution of Waring's problem for such groups, namely $w(G)^2 = G$ 
provided $|G|$ is large enough. Moreover, we prove a somewhat
stronger result, namely Theorem~\ref{bounded-waring} dealing
with two different words.

We begin with a basic geometric lemma:

\begin{prop}
\label{factor-through}
Consider the diagram
\begin{equation*}
\xymatrix{&X\ar[dl]_f\ar[dr]^g\\ Y\ar@{-->}[rr]&&Z \\}
\end{equation*}
of affine varieties over an algebraically closed field $k$.
If $Y$ is normal, $f$ is generically
smooth and surjective, and $g$ factors through $f$ at the level of closed points,
then $g$ factors through $f$ as a morphism of algebraic varieties.
\end{prop}

\begin{proof}
Let $A$, $B$, and $C$ denote the coordinate rings of $X$, $Y$, and $Z$ respectively.
They are Jacobson rings \cite[10.4.6]{EGA4-3}, so that every non-empty locally closed
subset of $X$, $Y$, or $Z$ contains at least one closed point \cite[10.1.2]{EGA4-3}.  Let
$K$ and $L$ denote the function fields of $X$ and $Y$ respectively.
We have the diagram
\begin{equation*}
\xymatrix{&L&A\ar@{_{(}->}[l]\\ 
K\ar@{^{(}->}[ur]&B\ar[ur]^<<<{\phi_f}\ar@{_{(}->}[l]&&C\ar@{-->}[ll]\ar[ul]_<<<{\phi_g} \\}
\end{equation*}
Surjectivity of $f$ implies that every prime ideal of $B$ 
contains $\ker\phi_f$, and as $B$ is an integral domain, this
means that $\phi_f$ is injective.
We regard $B$ as a subring of $A$ via $\phi_f$.  
Our goal is to prove that $\phi_g(C)\subset B$.  Suppose not.
Let $h\in A$ denote an element of $\phi_g(C)\setminus B$.
Thus $f$ factors through
$$f_h\colon \Spec B[h]\to\Spec B.$$
We consider four cases:

\begin{enumerate}
\item $L(h) = L$.
\item $L(h)$ is an inseparable extension of $L$.
\item $L(h)$ is a non-trivial separable extension of $L$.
\item $L(h)$ is a transcendental extension of $L$.
\end{enumerate} 

In case (1), as $B$ is normal and $h\in L\setminus B$, there exists a prime ideal $\p$ of $B$
of height $1$ such that $h\not\in B_{\p}$ \cite[Th.~38]{Ma}.  
In other words, $\p$ does not lie in the image of $f_h$,
contrary to the surjectivity hypothesis on $f$.

In case (2), the extension $L\to K$ is inseparable, contrary to the hypothesis that $f$ is 
generically smooth.

In case (3), the morphism $X\to \Spec B[h]$ is dominant (since $B[h]\subset A$).
By \cite[17.6.1]{EGA4-4},
there exists a dense open set $U_1\subset \Spec B[h]$ for which the morphism 
$f_h$ is
\'etale and by \cite[1.8.4]{EGA4-1}, a dense open set $U_2\subset \Spec B[h]$ contained in
the image of $X$ in $\Spec B[h]$.
Thus $\overline{f_h(U_1^c\cup U_2^c)}^c$ is a non-empty open set
in $\Spec B$ over which $f_h$ is \'etale
and every point in $\Spec B[h]$ lies in the image of $X$.
By \cite[15.5.9]{EGA4-3}, by passing to a smaller open set if necessary,
the cardinality of every fiber of $f_h$ equals $[L(h):L]$.
It follows that
there exist closed points $x_1,x_2\in X$ with the same image in $\Spec B$ but 
distinct images in $\Spec B[h]$.  If $\m_1$ and $\m_2$ are the corresponding maximal ideals
of $A$, then $h$ maps to different elements in $A/\m_1 = k$ and $A/\m_2 = k$.
Adding a suitable element of $k$ to $h$, we may assume that $h$ maps to $0$ in $A/\m_1$
but not in $A/\m_2$ and therefore some element of $C$ maps to $\m_1$ but not to $\m_2$.
It follows that $x_1$ and $x_2$ have distinct images in $\Spec C$, contrary to hypothesis.

The argument is the same in case (4), except that $\Spec B[h]\to \Spec B$ has the property that
the inverse image of every point in $\Spec B$ is infinite.

\end{proof}

\begin{prop}
\label{class-powers}
Let $G$ be a simply connected, almost simple algebraic group over an algebraically closed field $k$.
If $X\subset G$ is a conjugacy class of
positive dimension, then for all sufficiently large integers $n$, the morphism
$\pi_n\colon X^n\to G$ obtained by multiplying coordinates in $G$
is generically smooth and surjective.
\end{prop}

\begin{proof}
Let $X_i$ denote the Zariski closure of $\pi_i(X^i)$.  As $X^i$ is irreducible,
$X_i$ is the closure of a single point and therefore irreducible.  The sequence 
$\dim X_1,\dim X_2,\ldots$ is non-decreasing and must therefore stabilize, so we may
assume $\dim X_n = \dim X_{n+1} = \cdots$.   If $x_1,x_2\in X(k)$, then
$$x_1 X_n\cup x_2 X_n\subset X_{n+1},$$
so $x_1 X_n = x_2 X_n$, and 
\begin{equation}
\label{some-translations}
x_1^{-1} x_2 X_n = X_n.
\end{equation}

Let $Y_i$ denote the Zariski-closure of $\phi_i(X^{2i})$ (or, equivalently, the closure of
$\phi_i(X(k)^{2i})$, where
$$\phi_i(x_1,\ldots,x_{2i}) = x_1^{-1} x_2 x_3^{-1} x_4\cdots x_{2i-1}^{-1} x_{2i}.$$
If $n$ is sufficiently large, $\dim Y_m = \dim Y_{m+1}$, so
$$x_1^{-1} x_2 Y_m = x_1^{-1} x'_2 Y_m,$$
or
$$x_2^{-1} x'_2 Y_m = Y_m$$
for all $x_2,x'_2\in X(k)$.  This implies that $Y_m Y_m = Y_m$ and therefore $Y_m$ is
a connected algebraic subgroup of $G$.  As $X$ is conjugation-invariant, $Y_m$ is 
a normal subgroup,
and it follows that $Y_m = G$.  By (\ref{some-translations}),
$$X_n = Y_m X_n = G X_n = G.$$

Let $Z_i$ denote the Zariski closure of the complement of $\pi_i(X^i)$.
For $i\ge n$,
$Z_i$ is a proper closed subvariety of $G$.  If $x_1,x_2\in X(k)$, then
$$Z_{n+1}\subset x_1 Z_n \cap x_2 Z_n.$$
If $x_1 Z_n = x_2 Z_n$ for all $x_1,x_2\in X(k)$, then $Z_n$ is invariant under translation
by $Y_m = G$, which is impossible.
Thus
$x_1 Z_n\neq x_2 Z_n$ for some $x_1,x_2\in X(k)$.  If we form the vector of
dimensions of irreducible components of $Z_n$, arranged from largest to smallest, then
the vector of $Z_n$ majorizes that of $Z_{n+1}$.
We conclude that $Z_n=\emptyset$ for $n\gg 0$, and this gives surjectivity.

For generic smoothness, it suffices to find points $x_1,\ldots,x_n\in X(k)$ such that the map
$G^n\to G$ given by 
$$(g_1,\ldots,g_n)\mapsto g_1 x_1 g_1^{-1} g_2 x_2 g_2^{-1}\cdots g_n x_n g_n^{-1}$$
is smooth at the identity in $G^n$.  Equivalently, we must prove that
the space
\begin{multline*}
L(x_1,\ldots,x_n):= (1-\ad(x_1))(\g)+\ad(x_1)(1-\ad(x_2))(\g)+\cdots \\
+ \ad(x_1\cdots x_{n-1})(1-\ad(x_n))(\g)
\end{multline*}
equals $\g$ for suitable $x_1,\ldots,x_n$.  As 
$$L(x_1,\ldots,x_n) = (1-\ad(x_1))(\g) + \ad(x_1)(L(x_2,\ldots,x_n)),$$
we have
$$\dim L(x_1,\ldots,x_n)\ge \dim L(x_2,\ldots,x_n)$$
with equality if and only if 
$$(\ad(x_1^{-1})-1)(\g)\subset L(x_2,\ldots,x_n).$$

Now,
$$S_X:=\Span_{x\in X} (\ad(x^{-1})-1)(\g)$$
is invariant under the adjoint action of $G$.
As $x$ is not in the center of $G$, $S_X\neq 0$.  If the adjoint representation of
$\g$ is irreducible, which is usually the case since $G$ is almost simple, it follows that
$S_X = \g$.  To see that this is true in general, we observe that
since every element of $G$ is a product of elements of $X$,
$$S_X = \Span_{g\in G} (\ad(g) - 1)(\g),$$
and $G$ acts trivially on $\g/S_X$.  As $G$ is simply connected,
there is no non-trivial quotient of $\g$ on which $G$ acts trivially (see, e.g., \cite[1.11]{Pi}).

Therefore, the dimension determined by a sequence in $X$ can always be increased 
by prepending a suitable element of $X$ unless it already equals $\dim\g$. 
The proposition follows.
\end{proof}

We can now prove our basic irreducibility theorem: 

\begin{thm}
\label{basic-irred}
Let $w_1$ and $w_2$ denote non-trivial words in $n_1$ and $n_2$ letters respectively,
and let $w\in F_{n_1+n_2}$ denote their juxtaposition.
Let $G$ be a simply connected almost simple algebraic group over an
algebraically closed field $k$.
Then for all non-central elements $g\in G(k)$, the fiber $\pi_w^{-1}(g)$ is irreducible.
\end{thm}

We remark that the hypothesis that $G$ is simply connected is really needed.
Otherwise, we can take $w_1 = x_1 x_2 x_1^{-1} x_2^{-1}$ and
$w_2 = x_3 x_4 x_3^{-1} x_4^{-1}$ and the resulting morphism
$\pi_w$ factors through the universal cover of $G$.

\begin{proof}
We can express the fiber $\pi^{-1}_w(g)$ as the fiber product over $G$ of $\pi_{w_1}$
and a second morphism, namely the composition of $\pi_{w_2}$ with the 
``reflection'' map $x\mapsto gx^{-1}$.  We will prove that this fiber product is
geometrically irreducible.
 
Let $A$ denote the coordinate ring of $G$ and $K$ the field of fractions of $A$.
Let $K_i$ denote the separable closure of $K$ in the fraction field of
$A^{\otimes n_i}$, where the morphism $A\to A^{\otimes n_i}$ is that associated to
$\pi_{w_i}$.  Let $L_i$ be the splitting field of $K_i$ over $K$ and $\lambda_i\colon K\to L_i$
the natural inclusion map.

For any $g\in G(k)$, we define reflection, conjugation, and translation maps $G\to G$ as follows:
$$\rho_g(x) = gx^{-1},\ \chi_g(x) = gxg^{-1},\ \tau_g(x) = gx.$$
The induced automorphisms of $A$ and $K$ are also denoted $\rho_g$, $\chi_g$,
and $\tau_g$ respectively.  Let 
$$\chi_{g,i}(x_1,\ldots,x_{n_i}) = (gx_1 g^{-1},\ldots,gx_{n_i} g^{-1}).$$
We claim that if $g$ does not lie in the center of $G(k)$, then 
$\lambda_1$ and $\lambda_2\compose \rho_g$ give linearly disjoint finite extensions of
$K$.  As $A^{\otimes n_i}$ is geometrically irreducible over $K_i$ \cite[4.5.10]{EGA4-2}
and fiber products of geometrically irreducible schemes over a field are again
geometrically irreducible \cite[4.5.8]{EGA4-2}, this claim implies the theorem.
As $L_1$ and $L_2$ are Galois over $K$, it is equivalent to prove that
there do not exist subfields $\tilde K_1$ and $\tilde K_2$ of $L_1$ and $L_2$, properly containing
$K$, such that $\tilde K_1$ and $\tilde K_2$ are isomorphic as $K$-extensions.  

Suppose, on the contrary, that such an isomorphism 
$\iota\colon \tilde K_1\to\tilde K_2$ exists.  The injective homomorphisms
$\lambda_i$ factor through $\tilde K_i$, and we write $\kappa_i$ for the corresponding homomorphisms
$K\hookrightarrow \tilde K_i$.  Then we have the following diagram, in which the horizontal maps are isomorphisms and the vertical maps are finite field extensions:
\begin{equation*}
\xymatrix{\tilde K_1\ar[r]^{\chi_{h,1}}&\tilde K_1\ar[r]^\iota&\tilde K_2\ar[r]^{\chi_{h,2}^{-1}}
	&\tilde K_2\ar[r]^{\iota^{-1}}&\tilde K_1 \\
K\ar[u]_{\kappa_1}\ar[r]^{\chi_h}&K\ar[u]_{\kappa_1}\ar[r]^{\rho_g}&
	K\ar[u]_{\kappa_2}\ar[r]^{\chi_h^{-1}}&K\ar[u]_{\kappa_2}\ar[r]^{\rho_g^{-1}}&K\rlap{.}\ar[u]_{\kappa_1} \\}
\end{equation*}
The composition of horizontal arrows gives a diagram
\begin{equation*}
\xymatrix{\tilde K_1\ar[rr]^{\sigma_{g,h}}&&\tilde K_1 \\
K\ar[u]^{\kappa_1}\ar[rr]^{\tau_{hgh^{-1}g^{-1}}}&&K\rlap{,}\ar[u]_{\kappa_1} \\}
\end{equation*}
where $\sigma_{g,h}$ is an isomorphism depending on $g$ and $h$.  

Let $\tilde A$ denote the integral closure of $A$ in $\tilde K_1$.  
As $A$ is finitely generated over the field
$K$, it is a Japanese ring \cite[7.7.4]{EGA4-2}, so $\tilde A$ is a finitely generated $A$-module.  
The map $\kappa_1$ restricts to an injective ring homomorphism $A\to \tilde A$; by a slight abuse of
notation, we denote by $\kappa_1$ also the corresponding homomorphism of affine varieties
$\tilde G\to G$, where $\tilde  G=\Spec \tilde A$.
We fix $e\in \tilde G(k)$ lying over the identity in $G$.
For every $h_1,h_2\in G$, 
$\sigma_{g,h_1}$ and $\chi_{h_2,1}$ induce automorphisms of $\tilde G$.
Fix $h_1$ which does not commute with $g$ and let $h_2$ vary over $G$.
By taking the commutator of $\sigma_{g,h_1}$ and $\chi_{h_2,1}$ we obtain
a morphism $\Sigma\colon X\times \tilde G\to \tilde G$, where $X$ is the conjugacy class of 
$h_1 g h_1^{-1}g^{-1}$,
and each point of $X$ gives an automorphism of $\tilde G$ which covers the corresponding
translation of $G$.  We write $\Sigma_n$ for the morphism $X^n\times \tilde G\to \tilde G$ defined recursively by
$$\Sigma_n(x_1,\ldots,x_n,y) = \Sigma(x_1,\Sigma_{n-1}(x_2,\ldots,x_n,y)),$$
for $n\ge 2$ and $\Sigma_1 = \Sigma$.
Let $\Sigma_{n,e}\colon X^n\to \tilde G$ be defined by
$$\Sigma_{n,e}(x_1,\ldots,x_n) = \Sigma_n(x_1,\ldots,x_n,e).$$
The composition $\Sigma_{n,e}\compose\kappa_1$ gives the restriction to $X^n$
of the usual $n$-fold multiplication morphism $G^n\to G$.

By Proposition~\ref{class-powers}, the composition of $\Sigma_{n,e}\colon X^n\to \tilde G$
and $\tilde G\to G$ is surjective and generically smooth for all sufficiently large 
integers $n$.  It follows that $\Sigma_{n,e}$ itself is generically smooth and
dominant.  If the Zariski-closure $\tilde Z_n$ of the complement of
the image of $\Sigma_{n,e}$ is non-empty, then there exist $x_1,x_2\in X(k)$ such
that $\Sigma(x_1,\tilde Z_n)\neq \Sigma(x_2,\tilde Z_n)$ (because the images of these two subvarieties of $\tilde G$ in $G$ are distinct).  As 
$$\tilde Z_{n+1} \subset \Sigma(x_1,\tilde Z_n) \cap \Sigma(x_1,\tilde Z_n),$$
we see the vector of component dimensions of $\tilde Z_n$ majorizes that of
$\tilde Z_n$.  For $n$ sufficiently large, therefore, $\tilde Z_n$ is empty, i.e., 
$\Sigma_{n,e}$ is surjective.

We would like to apply Proposition~\ref{factor-through} to the upper triangle of
the diagram
\begin{equation}
\label{pentagon}
\xymatrix{&X^n\times \tilde G\ar[dl]_{\Sigma_{n,e}\times\id}\ar[dr]^{\Sigma_n} \\
\tilde G\times \tilde G\ar@{-->}[rr]^{\tilde\mu}\ar[d]_{\kappa_1\times \kappa_1}
	&&\tilde G\ar[d]^{\kappa_1} \\
G\times  G\ar[rr]&& G\rlap{.}\\}
\end{equation}
We observe first that the pentagon of solid arrows in (\ref{pentagon}) commutes because the
diagram
\begin{equation*}
\xymatrix{X\times \tilde G\ar[r]^<<<<<\Sigma\ar[d]_{\id\times \kappa_1}&\tilde G\ar[d]^{\kappa_1}\\ 
	X\times G\ar[r]& G\\}
\end{equation*}
commutes.  The variety $\tilde G\times \tilde G$ is normal since $\tilde G$ is so by construction.
To check that there exists a morphism $\tilde\mu$ at the level of closed points,
we note that if $x_1,\ldots,x_n$ belong to $X(k)$, the resulting map 
$\Sigma_n(x_1,\ldots,x_n,\cdot)$ is an automorphism of $\tilde G$ which covers
translation by $x_1\cdots x_n$.  Our claim asserts that, at least at the level of closed points,
this automorphism depends only on  $\Sigma_{n,e}(x_1,\ldots,x_n)$ and not on the actual
$n$-tuple $(x_1,\ldots,x_n)$.   Indeed, if we fix $\tilde g\in \tilde G$ and let $g = \kappa_1(\tilde g)$
there is a unique automorphism of $\tilde G$ covering $\tau_g$ and sending $e$ to $\tilde g$.
Thus, we may apply Proposition~\ref{factor-through} to define $\tilde\mu$.  
The resulting diagram commutes (i.e., the lower rectangle does so) because it commutes
at the level of closed points.

We claim that the arrow $\tilde \mu$ in (\ref{pentagon})
is a multiplication morphism making $(\tilde G,e)$ into an algebraic group
and $\tilde G\to G$ into a central isogeny.  As $G$ is simply connected, this implies that
$\kappa_1$ is an isomorphism \cite[1.5.4]{Ti}, and therefore that $K=\tilde K_1$, which proves the theorem.

To prove that $\tilde \mu$ satisfies the associativity axiom, we compare the two maps
$\tilde G^3\to \tilde G$ given by $\tilde \mu(\tilde \mu(x,y),z)$ and $\tilde \mu(x,\tilde \mu(y,z))$.  We have 
a diagram
\begin{equation*}
\xymatrix{\tilde G^3\ar@<2pt>[r]\ar@<-2pt>[r]\ar[d]_{\kappa_1^3}&\tilde G\ar[d]^{\kappa_1}\\
G^3\ar[r]&G\rlap{.}\\}
\end{equation*}
The closed subvariety of $\tilde G^3$ on which the two arrows agree is also open since
the fibers of $\tilde G\to G$ are finite.  It is non-empty since it contains $(e,e,e)$.
As $\tilde G^3$ is connected, $\tilde \mu$ is associative.

Next we construct the inverse on $\tilde G$.  We choose an automorphism of $\tilde G$
which covers the composition $\tau_{h^{-1}}\compose \rho_h$
(which maps $x\mapsto x^{-1}$).  This automorphism sends $e$ to some element of $\tilde G$ lying
over the identity of $G$.  By composing with a deck transformation, we can find a new automorphism
of $\tilde G$ which sends $e$ to itself and covers the inverse map on $G$.  Again we use connectedness to show that this morphism satisfies the diagram for inverse maps 
with respect to $\tilde \mu$.  Thus $(\tilde G,\tilde \mu,e)$ is an algebraic group and
the covering map $\tilde G\to G$ is a surjective homomorphism, i.e., an isogeny.  
As it is separable, it is central.  The theorem follows.

\end{proof}
 
We can now deduce Theorem~\ref{bounded-waring}

\begin{proof}
As in the proof of Proposition~\ref{fixed-rank}, it suffices to consider a fixed Dynkin diagram $\Delta$ 
and a group scheme $\pi\colon \cG\to\cS$ such that every simply connected,
almost simple algebraic group $G/\F_q$ with Dynkin diagram $\Delta$ is isomorphic to
$\cG_s$ for some $s\in \cS(\F_q)$.  Next we consider the morphism
$$\pi_w\colon \cG^{n_1+n_2}\to \cG$$
(where the fiber power is taken relative to $\cS$) given by a word $w$ which is a juxtaposition of
the words $w_1$ and $w_2$ taken in disjoint variables.  If $s\in \cS(\F_q)$ and $x\in \cG_s(\F_q)$
is a non-central element, then by Theorem~\ref{basic-irred}, 
the fiber of $\pi_w$ over $s$ is geometrically connected, so by Lemma~\ref{alg-geom}, if $q\gg 0$,
$x\in w(\cG_s(\F_q))$.  If $\Gamma$ denotes the quotient of $\cG_s(\F_q)$
by its center, it follows that the image of $x$ in $\Gamma$ lies in 
$w(\Gamma)$.
This accounts for all elements of $\Gamma$ except the identity, which 
trivially lies in $w(\Gamma)$.

\end{proof}

\section{Special word values in $\SL_2(\F_p)$}

In this section we take a closer look at 2-dimensional special
linear groups, which play a key role in our results for
alternating groups. Our main result below proves the existence
of elements of very high order in  $w(\SL_2(\F_p))$, under
certain conditions on $p$.

\begin{thm}
\label{nearly-max}
For every non-trivial word $w$ there exist constants $M_w$ and $N_w$ 
such that for every prime $p>N_w$ such that $p-1$ is not divisible
by $4$ nor by any prime $3\le \ell\le M_w$, there exists an element of order $\frac{p-1}2$
in $w(\SL_2(\F_p))$.
\end{thm}

\begin{proof}
Our main goal is to find $s\in\SL_2(\F_p)$ 
satisfying the following conditions:

\begin{itemize}

\item There exists $t\in\SL_2(\F_p)^d$ such that $\pi_w(t)=s$.

\item The centralizer of $s$ is a split torus.

\item There are two square roots of $s$ in $\SL_2(\F_p)$

\item If $k>1$ is not divisible by any prime $\le M_w$, then
$s$ has at most $k-1$ $k$th roots in 
$\SL_2(\F_p)$.

\end{itemize}

Let $P_k$ denote the $k$th Chebychev polynomial,
normalized so that
$$P_k(x+x^{-1}) = x^k+x^{-k}.$$
Note that $P_k(x)\in \Z[x]$.
If $\lambda\in\bar\F_p\setminus\{0\}$ and $k$ is an integer such that $\lambda^k\neq \lambda^{-k}$, then the
conditions $\{\lambda^k,\lambda+\lambda^{-1}\}\subset \F_p$ and $\lambda\in\F_p$ are equivalent.
Therefore, the specified conditions on $s$
are equivalent to the following conditions on the trace
$u=\tr(s)$:

\begin{itemize}

\item The equation $\tr(\pi_w(t))=u$ 
has at least one solution $t\in \SL_2(\F_p)^d$.

\item The equation $x+x^{-1} = u$ has two solutions in $\F_p$.

\item The equation $P_2(x)=u$ has two solutions in $\F_p$.

\item If $k>1$ is not divisible by any prime $\le M_w$,
	then the equation $P_k(x)=u$ has at most $k-1$ solutions
	$x\in \F_p$.

\end{itemize}

Let
$$\SL_2 = \Spec \Z[a_{11},a_{12},a_{21},a_{22}]/(a_{11}a_{22}-a_{12}a_{21}-1)$$
and $\A^1 = \Spec \Z[x]$.  
Let $\cX_1$ 
denote the integral closure of $\A^1$ in $\SL_2^d$, and write
$\psi\colon \SL_2^d\to\cX_1$ and $\phi_1\colon \cX_1\to \A^1$ for
the natural maps whose composition $\tau\colon \SL_2^d\to \A^1$
gives the trace of the morphism $\pi_w$ associated to $w$:
\begin{equation*}
\xymatrix{&\SL_2^d\ar[dl]_\psi\ar[dd]^\tau\ar[dr]^{\pi_w} \\
\cX_1\ar[dr]_{\phi_1}&&\SL_2\ar[dl]^{\tr} \\
&\A^1\rlap{.}}
\end{equation*}
The generic fiber
of $\psi$ is geometrically connected.  By the openness of geometric
connectivity \cite[9.7.7]{EGA4-3}, there exists a non-empty 
subscheme $\cU\subset\A^1$ 
such that 
$\psi|_{\tau^{-1}(\cU)}\to \phi_1^{-1}(\cU)$ has geometrically connected fibers.

Let 
$$\cX_2 = \Spec \Z[x,y]/(y^2-x^2+4),$$
and let $\phi_2\colon \cX_2\to \A^1$ denote the map $(x,y)\mapsto x$.  Let $\cX_3^{(k)} = \A^1$ for all
positive integers $k$, and $\phi_3^{(k)}\colon \cX_3^{(k)}\to\A^1$
be the Chebychev morphism given by $P_k$.  
The schemes $\cX_1$, $\cX_2$, and $\cX_3^{(k)}$ are generically
finite over $\A^1$, so shrinking $\cU$ if necessary we may assume
that over $\cU$, all three are finite \'etale.  Let
$M_w$ denote the degree of $\phi_1$.
If $N_w$ is sufficiently large, for every prime $p > N_w$, the 
$\F_p$-fiber of $\cU$ is non-empty.  For some $L_w$ independent of $p$,
for every $p>N_w$, $\F_p\setminus \cU(\F_p)$ has at most $L_w$ points.

We now fix $p$ sufficiently large and
let $U_1$, $U_2$, and $U_3^{(k)}$ denote the
inverse images of the $\F_p$-fiber $\cU_{\F_p}$
in $\cX_{1\F_p}$, $\cX_{2\F_p}$, and $\cX^{(k)}_{3\F_p}$
respectively.

We define the various fiber products of these curves as follows:
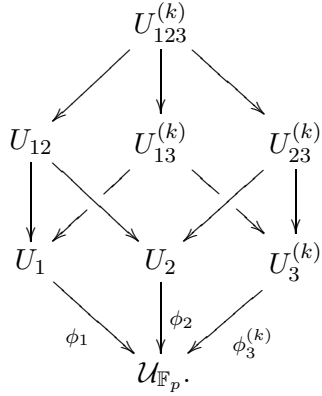
\begin{equation*}
\xymatrix{&U^{(k)}_{123}\ar[dl]\ar[d]\ar[dr] \\
U_{12}\ar[d]\ar[dr]&U^{(k)}_{13}\ar[dl]|\hole\ar[dr]|\hole
	&U^{(k)}_{23}\ar[dl]\ar[d] \\
U_1\ar[dr]_{\phi_1}&U_2\ar[d]^{\phi_2}&U^{(k)}_3\ar[dl]^{\phi_3^{(k)}} \\
&\cU_{\F_p}\rlap{.}}
\end{equation*}

We now rewrite our original set of conditions on $u\in\F_p$ (which
we now assume lies in the subset $\cU(\F_p)$):

\begin{itemize}
\item There exists $u_1\in U_1(\F_p)$ mapping to $u$.
\item There exists $u_2\in U_2(\F_p)$ mapping to $u$.
\item There exists $u_3^{(2)}\in U_3^{(2)}(\F_p)$ mapping to $u$.
\item If $k>1$ is not divisible by any prime $\le M_w$, the number of points
$u_3^{(k)}\in U^{(k)}_3(\F_p)$ mapping to $u$ is less than $k$.
\end{itemize}

The first condition may appear to be weaker than 
first of the previous set of conditions on $u$, but for 
$p\gg 0$, they are equivalent.  This is because
the fibers of $\psi$ over all points in $U_1(\F_p)$
are geometrically connected, so by Lemma~\ref{alg-geom}
they all have $\F_p$-points when $p\gg 0$.  The
second and third conditions mention only one preimage instead 
of two, but for two-to-one finite \'etale maps, there is no
difference.

We now suppose that $k$ is odd and $2k$ divides $p-1$.
By the Riemann-Hurwitz theorem, $U^{(2k)}_{123}$ is an open subset 
(the complement of at most $8kL_w$ points) of
a curve whose genus is $O(k)$.  (Here we use the fact that 
the \'etale morphism $U^{(2k)}_{123}\to U_{12}$ must extend to 
a tamely ramified cover of non-singular projective curves 
because its degree is less than $p$.)
By the Weil bound for curves, $U^{(2k)}_{123}(\F_p)$ has  
$p+O(k \sqrt{p})$ points.  As $2k$ divides $p-1$, 
a regular semisimple element in a split maximal torus 
of $\SL_2(\F_p)$ which admits one $2k$th root admits $2k$.
(In terms of our diagram, $U^{(2k)}_{123}$ is Galois over $U_{12}$.)
Therefore,
the image $I_k$ of $U^{(2k)}_{123}(\F_p)$ in
$U^{(2)}_{123}(\F_p)$
has $\frac{p}{k} + O(\sqrt{p})$ points.

The total number of points in $U^{(2)}_{123}(\F_p)$ is $p+O(\sqrt p)$.
Each maps to $u\in\cU(\F_p)$ which is the trace of 
an element of the form
$x^2$ in a split torus of $\SL_2(\F_p)$.  This element is the $r$th
power of a generator of the split torus for some $r$ dividing $p-1$.
We assume $p\equiv 3$ (mod $4$), 
so the highest power of $2$ dividing $r$ is $2$.
We would like to show that $x$ can be chosen so that $r=2$.
It suffices to show that for every odd divisor $k$ of $p-1$, $x$ is
not a $k$th power.  By hypothesis, $k$ has no prime divisors 
$\le M_w$, so every point on $U^{(2)}_{123}(\F_p)$ for which
the claim is false belongs to $I_k$.  Now,
$$\sum_{k|\frac{p-1}{2}} \mu(k) |I_k| = 2\phi\Bigl(\frac{p-1}{2}\Bigr) 
+ O\Bigl(\tau\Bigl(\frac{p-1}{2}\Bigr)\sqrt p\Bigr),$$
where $\mu$, $\phi$, and $\tau$ denote the M\"obius, Euler, and
number of divisors functions.  As $\tau(n) = o(n^{\epsilon})$ for
all $\epsilon > 0$, and $\phi(n) \ge n^{1-\epsilon}$ for all large
$n$, this implies that points of the desired form exist.
Let $u_{123}^{(2)}$ be one such, and let $u$ be its image
in $\cU(\F_p)$.  Then $u$ satisfies all the conditions 
needed to guarantee that it is the
trace of an element of $w(\SL_2(\F_p))$ of order $\frac{p-1}{2}$.

\end{proof}

\section{Alternating groups, I: three prime theorems}

In this section we prove Theorems \ref{cyc} and \ref{lower-bound}.
The proof combines our result on $\SL_2(\F_p)$ from the previous
section, with a certain approximation to the Goldbach conjecture
in the form of a three prime theorem.

To explain the idea, note that $\SL_2(\F_p)$ embeds in $A_{p+1}$,
and an element of order $(p-1)/2$ in $\SL_2(\F_p)$ has two non-trivial 
cycles and two fixed points in this embedding. Thus, if $n=p+1$
and $p$ satisfies the conditions of Theorem~\ref{nearly-max},
then we can produce $\sigma \in w(A_n)$ with $\cyc(\sigma) = 4$
(and with two fixed points).

To extend this idea for arbitrary $n$ we need to find suitable
primes $p_1, p_2, p_3$ such that $n = (p_1+1)+(p_2+1)+(p_3+1)+c$
where $c$ is a small constant, and then to embed
$\SL_2(\F_{p_1}) \times \SL_2(\F_{p_2}) \times \SL_2(\F_{p_3})$
in $A_n$, which produces a permutation with a bounded number
of cycles. Moreover, results from Additive Number Theory on 
the number of ways to write $n$ in the form above will provide us 
with information on the size of $w(A_n)$.  

We now give the detailed proof of Theorems \ref{cyc} and \ref{lower-bound}.

\begin{proof}
We are given a non-identity word $w$.
By Theorem~\ref{nearly-max}, for every prime $p > N_w$ such that
$p\equiv 3 \pmod{4}$ and $p-1$ not divisible by any odd prime $\le M_w$,
there exist $z_{p,1},\ldots,z_{p,d}\in \SL_2(\F_{p})$ such that
$w(z_{p,1},\ldots,z_{p,d})$ is of order $\frac{p-1}{2}$.
The action of $\SL_2(\F_p)$ on $\P^1(\F_p)$ determines an injective
homomorphism $\SL_2(\F_p)\to A_{p+1}$.  Let $\zeta_{p,i}$ denote the
image of $z_{p,i}$.  Thus $w(\zeta_{p,1},\ldots,\zeta_{p,d})$
consists of two $\frac{p-1}{2}$-cycles and two fixed points.

Next we claim that for every integer $N$
sufficiently large in terms of $M=M_w$, there exists
an integer $N'$ in the interval $[N-11,N]$ such that $N'-3$ is the sum of
three primes $p_i$ such that $p_i > N_w$, $p_i\equiv 3$ (mod $4$),
and $p_i-1$ is not divisible by any prime in $\le M_w$.
We do this by applying Ayoub's version of Vinogradov's three prime
theorem \cite{Ay}, asserting that every sufficiently large integer $N$ can be 
written in at least $c N^2 \log^{-3}N$ ways as a sum of three 
primes lying in any specified residue classes $a_1$, $a_2$, $a_3$ (mod $m$)
as long as $N$ is congruent to $a_1+a_2+a_3$ (mod $m$).  
For $N$ sufficiently large, this is enough to ensure that
$p_1$, $p_2$, $p_3$ can be chosen to be larger than $N_w$.

We choose $N'=12\lfloor N/12\rfloor$ so that $N'-3$ is the sum of
three primes each congruent to $11$ (mod $12$).  If $p_i$ is such a prime,
$p_i-1$ is not divisible by $3$.  For every prime $\ell$ 
greater than $3$, every residue class (mod $\ell$) can be written as the
sum of three residue classes from $2$ to $\ell-1$, so we can impose
(mod $\ell$) conditions on $p_1$, $p_2$, $p_3$ consistent with 
$\ell\nmid p_i-1$ and $p_1+p_2+p_3+3\equiv N'$ (mod $\ell$).
By the Chinese remainder theorem and Ayoub's result, 
the element 
$$w((\zeta_{p_1,1},\zeta_{p_2,1},\zeta_{p_3,1}),\ldots
(\zeta_{p_1,d},\zeta_{p_2,d},\zeta_{p_3,d})) \in 
A_{p_1+1}\times A_{p_2+1}\times A_{p_3+1}\subset A_{N'}$$
has at most $12$ cycles including $6$ fixed points, 
so by the inclusion $A_{N'}\subset A_N$,
there is an element in $w(A_N)$ with at most $23$ cycles,
of which at most $6$ are non-trivial.
We can construct a sequence of elements $\sigma_N\in A_N$ 
in this way for which the bound $M$ tends to $\infty$ as
$N$ goes to $\infty$.

This concludes the proof of Theorem~\ref{cyc}.

For Theorem~\ref{lower-bound}, we construct elements as above.
The centralizer of such an element in $A_n$ has order 
at most 
$$17! \frac{(p_1+1)^2(p_2+1)^2(p_3+1)^2}{16} = O(n^6)$$
Therefore, the conjugacy class of such an element has $O(n^{-6}|A_n|)$
elements.  As we have $O(n^{2-\epsilon})$ ways of choosing triples 
$(p_1,p_2,p_3)$, in all there are at least $O(n^{-4-\epsilon}|A_n|)$
elements in $w(A_n)$.

\end{proof}

\bigskip

\section{Alternating groups, II: a combinatorial construction}

In general it is difficult to say when a given conjugacy class
in a symmetric group appears in the product of two other
conjugacy classes.  In this section we present an explicit
construction that can be used to show that certain conjugacy
classes are contained in the squares of other classes.
This constitutes the elementary part of the proof of 
Theorem~\ref{squares}. 

We denote the number of fixed points of a permutation $\sigma \in S_n$
by $\fix(\sigma)$.

\begin{prop}
\label{naive-squares}
Let $n$ be a positive integer and $\alpha,\beta\in S_n$ permutations
belonging to conjugacy classes $A$ and $B$ respectively.
If 
$$\fix(\beta) \ge 7\cyc(\alpha),$$
and $\beta$ is even, then $B\subset A^2$.
\end{prop}

\begin{proof}
Let $[a,b]$ denote the sequence of consecutive integers
$a\le x\le b$; such a sequence will be called an \emph{interval
of length $1+b-a$}.  Given sequences 
$S_1,\ldots,S_r$, let $j(S_1,\ldots,S_r)$ denote their
concatenation.  
If $S$ is a sequence and $k$ is an integer, $T_k(S)$ is the sequence
obtained by adding $k$ to each element of $S$.
If $S$ is a $k$-term sequence of distinct integers
in $[1,n]$, let $(S)$ denote the corresponding $k$-cycle in $S_n$.
Let $0=a_0<a_1<\ldots<a_m=n$
be defined so that
$$\gamma = ([a_0+1,a_1])([a_1+1,a_2])\cdots([a_{m-1}+1,a_m])\in A.$$
We want to find $\delta\in A$ such that $\gamma\delta\in B$.

Let $b_k$ denote the number of $k$-cycles in $\beta$.  Let 
\begin{align*}
c_1&=b_1 = \fix(\beta)\ge 7\cyc(\alpha) = 7m,\\
c_2&=b_2,\\
c_3&=b_3+b_5+b_7+\cdots,\\
c_4&=b_4+b_6+b_8+\cdots,\\
d_3&=c_3,\\
d_4&=\lfloor c_2/2\rfloor,\\
d_6&=c_2 - 2d_4,\\
d_8&= \lfloor c_4/2\rfloor.
\end{align*}
As $b_2+b_4+b_6+\cdots$ is even, so is $c_2+c_4$, so 
$c_2$ and $c_4$ are both even or both odd and $c_4 = 2d_8+d_6$.

Our first task is a packing problem: to find a set $\Sigma$
consisting of $d_i$ intervals of
length $i$ for $i=3,4,6,8$,
each contained in some $[a_k+1,a_{k+1}]$, such that the resulting
intervals are mutually disjoint.  To show that
such a packing is possible, we iterate, starting from
$k=0$ and packing as many intervals of length $i\in\{3,4,6,8\}$ 
as possible (up to a limit of $d_i$ intervals of length $i$) into each
interval $[a_k+1,a_{k+1}]$ in turn, starting from
the right endpoint and working toward the left, so the union of intervals at
any point in the construction is ``right-justified.''  
For each $k$, we denote by $I_{k,1},I_{k,2},\ldots,I_{k,p_k}$
the successive elements of $\Sigma$ in $[a_k+1,a_{k+1}]$

If the packing process terminates before all the upper limits $d_i$
have been achieved, then the number of free spaces in all the
$[a_k+1,a_{k+1}]$ combined is $\le 7m$.
Since $b_1 = \fix(\beta)\ge 7m$, we have
$$3d_3+4d_4+6d_6+8d_8 = 2c_2+3c_3+4c_4\le b_2+b_3+b_4+\cdots \le n-7m,$$
so
$$|\{(k,l)\colon |I_{k,l}|=i\}| = d_i$$
for $i=3,4,6,8$.

We define $l_k$ so that for each $k$ from $0$ to $m-1$, the union of
the intervals $I_{k,1},I_{k,2},\ldots,I_{k,p_k}$ is $[l_k+1,a_{k+1}]$.
We write $e_{k,i}$ for the left endpoint of $I_{k,i}$
and define
\begin{equation}
\label{magic}
S_{k,i} = 
\begin{cases}
T_{e_{k,i}}(2,0,1)                       &|I_{k,i}|=3\\
T_{e_{k,i}}(3,0,1,2)                     &|I_{k,i}|=4\\
T_{e_{k,i}}(5,2,0,3,4,1)                 &|I_{k,i}|=6\\
T_{e_{k,i}}(7,5,2,0,3,6,4,1)             &|I_{k,i}|=8\\
\end{cases}
\end{equation}
Note in particular that the first term of $S_{k,i}$ is the last term of
$I_{k,i}$.
If $|I_{k,i}|=3$ (resp. $|I_{k,i}|=6$), we say that $e_{k,i}+1$ is a 
\emph{special point of odd type} (resp. \emph{even type}).
Intervals of length $4$ in $\Sigma$ have no special points.
If $|I_{k,i}|=8$, there are two special points of even type in the 
interval: $e_{k,i}+1$ and $e_{k,i}+6$.
We assign each special point a label, i.e., a positive integer $\lambda$,
in such a way that there are $b_{2\lambda+1}$ special points of odd type with 
label $\lambda$ and $b_{2\lambda+2}$ special points of even type 
with label $\lambda$.

Let $X_k$ denote the interval $[a_k+1,l_k]$ and $Y_k$ the subset of $X_k$
consisting of $a_k+x$ such that $x$ is odd and $a_k+x+1\le l_k$.
Thus $|X_k|\le 2|Y_k|+1$.  Setting 
$$X = X_0\cup X_1\cup \cdots \cup X_{m-1},\ 
Y = Y_0\cup Y_1\cup\cdots \cup Y_{m-1},$$
we have $|X| \le 2|Y|+m$.  Thus,
\begin{equation}
\label{enough}
\begin{split}
|Y| &\ge \frac{n-3d_3-4d_4-6d_6-8d_8-m}{2}
= \frac{n-2c_2-3c_3-4c_4 - m}{2} \\
&= \frac{-m+b_1}2+\sum_{i=5}^\infty \lfloor\frac{i-3}{2}\rfloor b_i 
 \ge 3m + \sum_{i=5}^\infty \lfloor\frac{i-3}{2}\rfloor b_i \\
\end{split}
\end{equation}

There exists a permutation $\epsilon\in S_n$
with the following properties:
\begin{itemize}
\item Every cycle of length $\ge 2$ contains exactly one element of
$[1,n]\setminus Y$, and this element is special.
\item The length of a cycle containing a special point of label $\lambda$ 
is $\lambda$.
\end{itemize}

The existence of $\epsilon$ follows from (\ref{enough}).

For $0\le k < m$, we define a cycle
$$\delta_k = (j(S_{k,1},\ldots,S_{k,p_k},R_k)),$$
where $R_k$ is the reverse of the sequence $[a_k+1,l_k]$.
Thus, $\delta_0\delta_1\cdots\delta_{m-1}$ belongs to $A$.
We define
$$\delta = \epsilon \delta_0\cdots\delta_{m-1}\epsilon^{-1}\in A.$$
We claim $\gamma\delta$ belongs to $B$, i.e., that 
it has exactly $b_j$ $j$-cycles for each 
$j$.  As $\sum_j j b_j = n$, it suffices to prove that there
are at least $b_j$ $j$-cycles for each $j$.
There are four cases to check: $j=1$, $j=2$, $j\ge 3$ odd, and $j\ge 4$
even.

For $j=1$, we note that every element of $x\in R_k$ such that
$x$ and $x-1$ are fixed by $\epsilon$ satisfies
$\gamma\delta(x) = \gamma(x-1) = x$. Thus,
\begin{equation*}
\begin{split}
\fix(\gamma\delta) 
&\ge n - \sum_{k,i} |I_{k,i}| - 2|Y\setminus \fix(\epsilon)|\\
& = n - 3d_3 - 4d_4 - 6d_6 - 8d_8 -2|Y\setminus \fix(\epsilon)|\\
&= n - 2c_2 - 3c_3 - 4c_4 - 2|Y\setminus \fix(\epsilon)|\\
&= n - 2b_2 - 3b_3 - 4b_4 - 3b_5 - 4b_6 - \cdots - 2|Y\setminus\fix(\epsilon)|\\
&= n - 2b_2 - 3b_3 - 4b_4 - 5b_5 - 6b_6 - \cdots \\
&= b_1.
\end{split}
\end{equation*}

For $j=2$ we consider all $I_{k,i}\in \Sigma$ which are of length $4$.
For each such interval, $\gamma\delta$ maps
\begin{align*}
e_{k,i}&\mapsto\gamma\delta_k(e_{k,i})=\gamma(e_{k,i}+1)=e_{k,i}+2\\
e_{k,i}+1&\mapsto\gamma\delta_k(e_{k,i}+1)=\gamma(e_{k,i}+2)=e_{k,i}+3\\
e_{k,i}+2&\mapsto\gamma\delta_k(e_{k,i}+2)=\gamma(e_{k,i}-1)=e_{k,i}\\
e_{k,i}+3&\mapsto\gamma\delta_k(e_{k,i}+3)=\gamma(e_{k,i})=e_{k,i}+1.\\
\end{align*}
This produces a total of $2d_4$ $2$-cycles.  If there is an interval
$I_{k,i}$ of length $6$, then $\gamma\delta$ maps
\begin{align*}
e_{k,i}+5&\mapsto\gamma\delta_k(e_{k,i}+5)=\gamma(e_{k,i}+2)=e_{k,i}+3\\
e_{k,i}+3&\mapsto\gamma\delta_k(e_{k,i}+3)=\gamma(e_{k,i}+4)=e_{k,i}+5.\\
\end{align*}
Thus, $\gamma\delta$ has at least $2d_4+d_6=c_2=b_2$ cycles of length $2$.

If $j\ge 3$ is odd, let $\lambda = (j-1)/2$.  There are $b_j$
special points of odd type with label $\lambda$.  Each belongs to
an interval $I_{k,i}$ of length $3$.  The special point in such an interval
is $e_{k,i}+1$.  Thus $\gamma\delta$ maps
\begin{align*}
e_{k,i}&\mapsto\epsilon(e_{k,i}+1)+1\\
\epsilon^r(e_{k,i}+1)+1&\mapsto\epsilon^{r+1}(e_{k,i}+1)+1,\ 1\le r<\lambda\\
e_{k,i}+2&\mapsto e_{k,i}+1\\
\epsilon^{r+1}(e_{k,i}+1)&\mapsto\epsilon^r(e_{k,i}+1),\ 1\le r < \lambda \\
\epsilon(e_{k,i}+1)&\mapsto e_{k,i}\\
\end{align*}
producing a $j$-cycle.  It follows that $\gamma\delta$ has at least $b_j$ $j$-cycles.

If $j\ge 4$ is even, let $\lambda = j/2-1$.  There are
$b_j$ special points of even type with label $\lambda$.
Each special point with label $\lambda$ in an interval $I_{k,i}$ of length $8$
contributes a $j$-cycle because $\gamma\delta$ maps
\begin{align*}
e_{k,i}&\mapsto e_{k,i}+4\\
e_{k,i}+4&\mapsto\epsilon(e_{k,i}+1)+1\\
\epsilon^r(e_{k,i}+1)+1&\mapsto\epsilon^{r+1}(e_{k,i}+1)+1,\ 1\le r<\lambda\\
e_{k,i}+2&\mapsto e_{k,i}+1\\
\epsilon^{r+1}(e_{k,i}+1)&\mapsto\epsilon^r(e_{k,i}+1),\ 1\le r < \lambda \\
\epsilon(e_{k,i}+1)&\mapsto e_{k,i}\\
\end{align*}
and
\begin{align*}
e_{k,i}+3&\mapsto\epsilon(e_{k,i}+6)+1\\
\epsilon^r(e_{k,i}+6)+1&\mapsto\epsilon^{r+1}(e_{k,i}+6)+1,\ 1\le r<\lambda\\
e_{k,i}+7&\mapsto e_{k,i}+6\\
\epsilon^{r+1}(e_{k,i}+6)&\mapsto\epsilon^r(e_{k,i}+6),\ 1\le r < \lambda \\
\epsilon(e_{k,i}+6)&\mapsto e_{k,i}+5\\
e_{k,i}+5&\mapsto e_{k,i}+3.\\
\end{align*}
Likewise, if $I_{k,i}$, of length $6$, has a special point with label $\lambda$,
then $\gamma\delta$ maps
\begin{align*}
e_{k,i}&\mapsto e_{k,i}+4\\
e_{k,i}+4&\mapsto\epsilon(e_{k,i}+1)+1\\
\epsilon^r(e_{k,i}+1)+1&\mapsto\epsilon^{r+1}(e_{k,i}+1)+1,\ 1\le r<\lambda\\
e_{k,i}+2&\mapsto e_{k,i}+1\\
\epsilon^{r+1}(e_{k,i}+1)&\mapsto\epsilon^r(e_{k,i}+1),\ 1\le r < \lambda \\
\epsilon(e_{k,i}+1)&\mapsto e_{k,i}.\\
\end{align*}
Thus, again we have at least $b_j$ $j$-cycles, and this finishes the proof.
\end{proof}

\bigskip
\section{Alternating groups, III: characters and probability}

Here we prove a character bound for symmetric groups which
may have some independent interest, and combine it with 
Proposition~\ref{naive-squares} and
a probabilistic argument to prove 
that conjugacy classes of permutations
which contain a small number of cycles 
have the property that their squares cover
all of $A_n$.  

This proves Theorem~\ref{squares} and
answers an old question of J.~L.~Brenner
(but falls short of his strong conjecture \cite{Br}).
For our purposes, however, it is more than enough;
combining it with Theorem~\ref{cyc}, we deduce Theorem~\ref{two-words},
giving a best possible solution for Waring's problem for $A_n$.
This will complete the proof of all results stated in the
Introduction.

Recall that the irreducible characters $\chi \in Irr(S_n)$
are parameterized by partitions $\la$ of $n$, which correspond
to Young diagram. Let $\chi_{\la}$ be the character corresponding
to a diagram $\la$. 

We define {\it layers} of a diagram $\la$ as follows. 
The first layer consists of the first row and column of $\la$.
Once we remove the first layer we obtain a smaller diagram
$\la '$. The layers of $\la$ are defined recursively as
the layers of $\la '$ together with the first layer of $\la$.

We shall use the Murnaghan-Nakayama Rule \cite[21.1]{Ja}.
By a {\it rim $r$-hook} $\nu$ in a $\la$-diagram, we mean a connected part of
the rim containing $r$ nodes, which can be removed to leave a proper diagram,
denoted by $\la \backslash \nu$. 
If, moving from right to left, the rim hook $\nu$ 
starts in row $i$ and finishes in column $j$, then the {\it leg-length}
$l(\nu)$ is defined to be the number of nodes below the $ij$-node
in the $\la$-diagram.

\begin{prop}\label{murn} {\rm (Murnaghan-Nakayama Rule)} 
Let $\rho \pi \in S_n$, where
$\rho$ is an $r$-cycle and $\pi$ is a permutation of the remaining $n-r$
points. Then 
\[
\chi_\la (\rho \pi) = \sum_\nu (-1)^{l(\nu)} \chi_{\la \backslash \nu} (\pi),
\]
where the sum is over all rim $r$-hooks $\nu$ in a $\la$-diagram.
\end{prop}

We shall now bound the character values of a permutation as
a function of its number of cycles alone.

\begin{thm} 
\label{char-bound}
Let $\sigma \in S_n$ be a permutation with $k$
cycles (including $1$-cycles). Then 
\[
|\chi(\sigma)| \le 2^{k-1}k!
\]
for all irreducible characters $\chi$ of $S_n$.
\end{thm}

\begin{proof}

We argue by induction on $k$. For $k=1$ our permutation 
$\sigma$ is an $n$-cycle, and it is a well known 
consequence of Proposition~\ref{murn} that
$|\chi(\sigma)| \le 1$, and moreover, if $\chi(\sigma) \ne 0$
then the Young diagram of $\chi$ consists of a single
layer.

Suppose now that $k \ge 2$ and let $\chi = \chi_{\la}$.
Write $\sigma = \rho \pi$ where $\rho$ is an $r$-cycle
and $\pi$ is a permutation on the remaining $n-r$ points.
We now apply the Murnaghan-Nakayama Rule above.
Since $\pi$ decomposes into $k-1$ cycles induction yields
\[
|\chi_{\la \backslash \nu} (\pi)| \le 2^{k-2}(k-1)!
\]
for all rim $r$-hooks $\nu$ in the diagram $\la$.

To complete the proof we may assume $\chi_{\la}(\sigma) \ne 0$.

\medskip

\noindent {\bf Claim:} The number rim $r$-hooks in $\la$ is at most $2k$. \par
\smallskip
To show this first note that, by repeatedly applying the
Murnaghan-Nakayama Rule, we see that $\la$ has at most $k$ layers
(otherwise $\chi_{\la}(\sigma) = 0$).

Now, each rim $r$-hook $\nu$ of $\la$ has a starting point, its
rightmost upmost point, and an end point, its leftmost downmost
point. Each of these points determines the rim $r$-hook uniquely
(by going $r$ steps on the boundary in the suitable direction).
We define the starting point and the end point of each layer of 
the diagram $\la$ in a similar manner.

Since the remaining part $\la \backslash \nu$ of the diagram should 
be connected, either the starting point of $\nu$ is the starting
point of some layer of $\la$, or the end point of $\nu$ is the end 
point of some layer of $\la$. Since there are at most $k$ layers, it now
follows that the number of rim $r$-hooks in $\la$ is at most
$2k$, proving the claim.

\medskip

The Murnaghan-Nakayama Rule now expresses $\chi_{\la}(\sigma)$
as a sum of at most $2k$ terms, each of absolute value at most
$2^{k-2}(k-1)!$. This yields
\[
|\chi_{\la}(\sigma)| \le 2k \cdot 2^{k-2}(k-1)! = 2^{k-1}k!,
\]
completing the proof of the theorem.

\end{proof}

Note that $|\chi(\sigma)|$ may be close to $(k!)^{1/2}$,
for instance take $\sigma = 1$ (so $k = n$) and $\chi$ a 
character of highest degree.
This shows that the upper bound in Theorem~\ref{char-bound} cannot be
significantly improved, at least for large $k$.

Given conjugacy classes $C_i = x_i^G$ ($i = 1,2$)
in a finite group $G$, let $P_{C_1, C_2}$ denote 
the distribution of
the random variable $y = y_1y_2$, where $y_i \in C_i$
is randomly chosen with uniform distribution on $C_i$. 
The following is standard.

\begin{lem}
\label{convolution}
With the above notation we have, for $g \in G$,
\[
P_{C_1,C_2}(g) = |G|^{-1} \sum_{\chi \in Irr(G)} 
\frac{\chi(x_1)\chi(x_2)\chi(g^{-1})}{\chi(1)}.
\] 
\end{lem}

For a finite group $G$ define
\[
\zeta_G(s) = \sum_{\chi \in Irr(G)} \chi(1)^{-s}
\]
where $s$ is a real number.

We need the following two results from \cite{LiSh2},
related to character degrees.

\begin{lem} 
\label{Witten}
For $s>0$ we have $\zeta_{S_n}(s) = 2 + O(n^{-s})$.
\end{lem}

\begin{lem}
\label{dim-bound}
Let $\lambda = \lambda_1, \ldots , \lambda_m$ be a partition of $n$ 
such that $\lambda_1 \ge m$.
Let $t = n-\lambda_1$.  Then

(i) $\chi_\lambda(1) \ge {{n-t} \choose t}$.

(ii) If $t \ge \epsilon n$ for some $\epsilon > 0$ then   
$\chi_\lambda(1) \ge c^n$ where $c>1$ depends on $\epsilon$. 

\end{lem}

We can now prove Theorem~\ref{squares}.

\begin{proof}
Throughout this proof, we assume without further comment that $n$
is sufficiently large.
Let $\sigma$ be a permutation in $S_n$ with 
$$\cyc(\sigma) \le n^{1/128}.$$
We have to show that $(\sigma^{S_n})^2 = A_n$.

We first claim that
\begin{equation}
\label{sigma-estimate}
|\chi(\sigma)| \le \chi(1)^{1/63},
\end{equation}
for all irreducible characters $\chi$ of $S_n$.

This is clear if $\chi(1)=1$ so suppose $\chi = \chi_{\lambda}$
is a non-linear character, where $\lambda = \lambda_1, \ldots , \lambda_m$
is a partition. 
Since we are only interested in absolute values, we may
replace a partition by its transpose and assume that 
$\lambda_1 \ge m$ (so Lemma~\ref{dim-bound} is applicable).

Write $\lambda_1 = n - t$ and $k = cyc(\sigma)$. Then
Theorem~\ref{char-bound} and our assumption on $k$ yield
\begin{equation}
\label{char-val}
|\chi_\lambda(\sigma)| \le 2^{k-1}k! \le (n^{1/128})^{n^{1/128}}.
\end{equation}

If $t \ge n/3$, Lemma~\ref{dim-bound}(ii) gives
an exponential lower bound on $\chi_\lambda(1)$,
so (\ref{sigma-estimate}) follows.

If $n^{1/127} \le t < n/3$, then Lemma~\ref{dim-bound}(i) yields
\[
\chi(1) \ge {{n-t} \choose t} \ge 2^t \ge 2^{n^{1/127}}.
\]
Combining this with (\ref{char-val}) we obtain (\ref{sigma-estimate}).

Finally, suppose $1\le t < n^{1/127}$. Let $f$ denote the number of
fixed points of $\sigma$.
We use \cite[(14)]{MS} to show that
$$|\chi(\sigma)| \le t! \sum_{a,b \ge 0, a+2b \le t}  
{f \choose a}{{k-f} \choose b} \le t! \sum_{0 \le i \le t} {k \choose i}
\le t! \cdot (t+1) \cdot k^t.$$

Using the bounds $k \le k_0:=n^{1/128}$ and $t \le t_0:=n^{1/127}$ this
gives
\[
|\chi(\sigma)| \le (t_0 k_0)^t \le n^{(1/128+1/127)t}.
\]
It is easy to see that for $t \le n^{1/127}$ and fixed $\delta > 0$ we have
\[
 n^{(1-1/127 -\delta)t} \le {{n-t} \choose t}.
\]
Choosing $\delta$ such that $1/128 + 1/127 = (1-1/127 -\delta) \cdot 1/63$ 
and noting that $\delta > 0$ it follows that
\[
|\chi(\sigma)| \le {{n-t} \choose t}^{1/63} \le \chi(1)^{1/63}.
\]
This completes the proof of (\ref{sigma-estimate}).

Now, let $C=\sigma^{S_n}$ and let $\pi$ denote an even permutation.

If $\pi$ has more than $7n^{1/128}$ fixed points,
then we have $\pi\in C^2$ by Proposition~\ref{naive-squares}.
We therefore assume from now on that $\fix(\pi) < 7n^{1/128}$.
By \cite[Theorem~B]{MS}, this implies
$$|\chi(\pi^{-1})|\le |\chi(1)|^{30/31}.$$

By Lemma~\ref{convolution} we have
\[
P_{C,C}(\pi) = \frac{1}{n!} \sum_{\chi \in Irr(S_n)} 
\frac{{\chi(\sigma)^2\chi(\pi^{-1})}} 
{\chi(1)}.
\]
The condition on the sign of $\pi$ shows that the contribution
of the linear characters of $S_n$ to $P_{C,C}(\pi)$
is $2/n!$.   On the other hand,
$$\Bigm|\sum_{\{\chi\colon \chi(1)>1\}} 
\frac{{\chi(\sigma)^2\chi(\pi^{-1})}}
{{\chi(1)}}\Bigm|
\le  \sum_{\{\chi\colon \chi(1)>1\}} 
\frac{{(\chi(1)^{1/63})^2 \chi(1)^{30/31}}} 
{{\chi(1)}} =  \sum_{\{\chi\colon \chi(1)>1\}} \chi(1)^{-s},
$$
where $s = 1/31 - 2/63 > 0$. 
Combining this with Lemma~\ref{Witten} we conclude that
\[
|P_{C,C}(\pi) - 2/n!| \le \frac 1{n!} (\zeta_{S_n}(s)-2)
= O(n^{-s}/n!).
\]
This implies $P_{C,C}(\pi) > 0$, so $\pi \in C^2$.

The theorem is proved.

\end{proof}

Combining various tools above,
we can finally prove Theorem~\ref{two-words}.

\begin{proof}

We have to show, given words $w_1, w_2 \ne 1$, that 
$w_1(A_n)w_2(A_n) = A_n$ for all $n \ge N(w_1, w_2)$.

We use Theorem~\ref{cyc} and its notation.
Let $N = N(w_1, w_2) := \max_{i=1,2} N(w_i)$, and
let $n \ge N$. Then Theorem~\ref{cyc} shows that 
there exists a permutation $\sigma_n \in w_1(A_n) \cap w_2(A_n)$ 
consisting of at most 23 cycles (of which at least $6$ are of length 1. 
so the $A_n$-conjugacy class of $\sigma_n$ is the same as its 
$S_n$-conjugacy class). 
Clearly $\sigma_n^{S_n} \subseteq  w_1(A_n) \cap w_2(A_n)$
and by Theorem~\ref{total-squares} we have
\[
w_1(A_n)w_2(A_n) \supseteq (w_1(A_n) \cap w_2(A_n))^2
\supseteq (\sigma_n^{S_n})^2 = A_n.
\]
The theorem is proved.

\end{proof}

\bigskip

\section{Intersection theorems}

Our methods enable us to deduce somewhat surprising results
concerning intersections of the form $\cap_{i=1}^k w_i(G)$,
where $w_1, \ldots , w_k$ are given words. 
Note that in free groups such intersection may well be
trivial. However, we can prove that in the groups considered
in this paper such intersections are very large.
More precisley, we have

\begin{thm} Let $k \ge 1$ and let $w_1, \ldots , w_k$ be
non-trivial words.

(i) There exists $N = N(w_1, \ldots , w_k)$ such that for all
$n \ge N$ we have
\[
(w_1(A_n) \cap \ldots \cap w_k(A_n))^2 = A_n;
\]

(ii) For every $\epsilon > 0$ there exists 
$N = N(w_1, \ldots , w_k, \epsilon)$ such that
\[
|w_1(A_n) \cap \ldots \cap w_k(A_n)| \ge n^{-4-\epsilon} |A_n|.
\]
\end{thm}

Indeed, part (i) follows from Theorems~\ref{squares} and~\ref{cyc},
using the fact that, for $n \gg 0$, 
$\sigma_n^{S_n} \subseteq \cap_{i=1}^k w_i(A_n)$.

The proof of part (ii) is similar to that of Theorem~\ref{lower-bound}.

Our next result is an intersection theorem for groups of Lie type.

\begin{thm} Let $k \ge 1$ and let $w_1, \ldots , w_k$ be
non-trivial words.
Let $G$ a simply connected almost simple 
algebraic group
over a finite field $\F_q$.  

(i) There is a positive constant $c$ depending
only on $w_1, \ldots , w_k$ and $\dim G$ such that 
\[
|w_1(G(\F_q)) \cap \ldots \cap w_k(G(\F_q))| \ge c|G(\F_q)|.
\]

(ii) There is a number $N$ depending only on
$w_1, \ldots , w_k$ such that, if $G$ has rank $r$ and is not of 
type $A_r$ or $^2 A_r$, and $|G(\F_q)| \ge N$, then 
\[
|w_1(G(\F_q)) \cap \ldots \cap w_k(G(\F_q))| \ge cr^{-1}|G(\F_q)|,
\]
where $c>0$ is an absolute constant.

(iii) There is a number $N$ depending only on
$w_1, \ldots , w_k$ and $\dim G$ such that, if 
$\Gamma = G(\F_q)/Z(G(\F_q))$ is a finite simple group
associated with $G$, and $|\Gamma| \ge N$, then
\[
(w_1(\Gamma) \cap \ldots \cap w_k(\Gamma))^2 = \Gamma.
\]

\end{thm}

The proofs of all three statements are modelled on the $k=1$ case, replacing the morphism
$w\colon G^n\to G$ (or in the case of (ii), $w\colon \SL_2^n\to\SL_2$) with the fiber product, relative to $G$,  of  the morphisms $w_i\colon G^{n_i}\to G$ (respectively $w_i\colon \SL_2^{n_i}\to\SL_2$).

\end{document}